\documentclass[12pt]{article}
\usepackage{amssymb,amsmath,mathrsfs,amsthm,latexsym}
\usepackage{hyperref}
\usepackage{graphicx}
\usepackage{color}
\usepackage[OT2,OT1]{fontenc}
\usepackage{upquote}
\usepackage[numbers,sort&compress]{natbib}
\usepackage{alphabeta}
\usepackage{alltt} 

\usepackage{xcolor}
\usepackage{framed}
\definecolor{shadecolor}{cmyk}{0,0,0,0.05} 

\newtheorem{theorem}{T{\hskip 0pt\footnotesize\bf HEOREM}}[section]
\newtheorem{lemma}[theorem]{L{\hskip 0pt\footnotesize\bf EMMA}}

\begin{document}

\title{\Large\bf A proof of the irrationality of $\pi$ motivated by nested radicals with roots of $2$}

\author{
\normalsize\bf Sanjar M. Abrarov, Rehan Siddiqui, Rajinder Kumar Jagpal \\
\normalsize\bf and Brendan M. Quine}

\date{September 23, 2026}
\maketitle

\begin{abstract}
In this work, we prove the irrationality of $\pi$ motivated by nested radicals with roots of $2$ of the kind $c_k = \sqrt{2 + c_{k - 1}}$ and $c_0 = 0$. Sample computations showing how the rational approximation tends to $\pi$ as the integer $k$ increases are presented.
\vspace{0.2cm}
\\
\end{abstract}

\section{Introduction}

In 1714, the English mathematician Roger Cotes discovered a remarkable identity \cite{Stillwell2010, Wilson2018}
\[
ix = \ln\left(\cos(x) + i\sin(x)\right).
\]
A few decades later, the Swiss mathematician Leonhard Euler independently rediscovered a reformulated form of this identity as
\[
e^{ix} = \cos(x) + i\sin(x)
\]
from which it follows that
\[
e^{i\pi} + 1 = 0.
\]
This equation, also known as Euler's identity, is commonly considered the most beautiful formula in mathematics as it relates the ubiquitous constants $\pi$ and $e$ to each other \cite{Wilson2018}. Sometimes these constants $\pi$ and $e$ are also regarded as Archimedes' constant and Euler's number, respectively.

A proof of the irrationality of the constant $e$ may not be difficult (see for example \cite{Coolidge1950, Davidson2023}). However, it was not easy to find a proof of the irrationality of $\pi$; a long time passed since the discovery of $\pi$ by ancient Babylonians and Egyptians \cite{Beckmann1971, Berggren2004, Agarwal2013} before its irrationality was proven.

The first proof that $\pi$ is irrational was given by the Swiss mathematician Johann Heinrich Lambert in 1761 \cite{Berggren2004, Angell2022} (see also \cite{Laczkovich1997}). In his work, Lambert showed that if $x \ne 0$ in the following infinite continuous fraction
\begin{equation}
\label{CF4TF} 
\tan(x) = \cfrac{x}{1 - \cfrac{x^2}{3 - \cfrac{x^2}{5 - \cfrac{x^2}{7 - \cfrac{x^2}{9 -  \ddots}}}}},
\end{equation}
then the value of $x$ cannot be rational when the expansion on the right side is rational. Therefore, in the equation
\[
\tan\left(\frac{\pi}{4}\right) = 1
\]
the constant $\pi$ must be irrational.

The first proof of the irrationality of $\pi$ by contradiction was found in 1873 by the French mathematician Charles Hermite \cite{Zhou2011}. There are several other proofs of the irrationality of $\pi$ \cite{Jeffreys2011, Niven1947, Huylebrouck2001, Bourbaki2004, Roegel2020, Chow2024}. One of them, published by Niven in 1947, is particularly interesting and attracts much attention. In his work \cite{Niven1947}, Niven proved the irrationality of $\pi$ also by contradiction. In particular, with the help of the series expansion
\[
F(x) = \sum_{m = 0}^n (-1)^m \frac{d^{\,2m}}{dx^{2m}} f(x),
\]
where
\[
f(x) = \frac{x^n(a - bx)^n}{n!},
\]
he showed that it is impossible to represent $\pi$ as a ratio of two integers $a$ and $b$. Despite a long history, research on the irrationality of $\pi$ still remains interesting \cite{Damini2020, Roegel2020, Angell2022, Chow2024}.

In this work, we present a proof of the irrationality of $\pi$ motivated by nested radicals of the kind $c_k = \sqrt{2 + c_{k - 1}}$, where $c_0 = 0$. Nested radicals of this kind have been used in our earlier publications \cite{Abrarov2022, Abrarov2023} to generate Machin-like formulas for $\pi$. To the best of our knowledge, this proof is new and has not been previously reported.

The outline of the remaining parts of this article is as follows. Section 2 presents preliminaries, Section 3 shows the motivation and a proof of the irrationality of $\pi$, and Section 4 includes examples of rational approximations.

\section{Preliminaries}

The identity \eqref{AF4P} below has been used in our previous publications \cite{Abrarov2021, Abrarov2022, Abrarov2023} as a starting point to generate Machin-like formulas that can be algorithmically implemented for rapid and efficient computation of the decimal digits of $\pi$ \cite{Agarwal2013}. The following theorem \ref{T4AF} shows how this identity can be derived.

\begin{theorem}
\label{T4AF} 
The following equation \cite{Abrarov2018}
\begin{equation}
\label{AF4P} 
\frac{\pi}{4} = 2^{k - 1}\arctan\left(\frac{\sqrt{2 - c_{k - 1}}}{c_k}\right), \qquad k \ge 1,
\end{equation}
holds.
\end{theorem}

\begin{proof}
Using the double-angle identity
\[
\cos(2x) = 2\cos^2(x) - 1,
\]
by induction it follows that
\[
\cos\left(\frac{\pi}{2^2}\right) = \frac{1}{2}{\sqrt 2} = \frac{1}{2}{c_1},
\]
\[
\cos\left(\frac{\pi}{2^3}\right) = \frac{1}{2}\sqrt{2 + \sqrt 2} = \frac{1}{2}{c_2},
\]
\begin{equation}
\label{CI} 
\cos\left(\frac{\pi }{2^{k + 1}}\right) = \frac{1}{2}\underbrace{\sqrt {2 + \sqrt{2 + \sqrt{2 + \ldots  + \sqrt 2}}}}_{k \,\, \text{square roots}} = \frac{1}{2}{c_k}.
\end{equation}
Therefore, we get
\begin{equation}
\label{SI} 
\sin\left(\frac{\pi}{2^{k + 1}}\right) = \sqrt{1 - \cos^2\left(\frac{\pi}{2^{k + 1}}\right)} = \sqrt{1 - \frac{1}{4}c_k^2}.
\end{equation}

Thus, using equations \eqref{CI} and \eqref{SI}, we obtain
\[
\begin{aligned}
\tan\left(\frac{\pi}{2^{k + 1}}\right) &= \frac{\sqrt{1 - \cos^2\left(\frac{\pi}{2^{k + 1}}\right)}}{\cos\left(\frac{\pi}{2^{k + 1}}\right)} \\
&= \frac{\sqrt{1 - \frac{1}{4}c_k^2}}{\frac{1}{2}c_k} = \frac{\sqrt{2 - c_{k - 1}}}{c_k}
\end{aligned}
\]
or
\[
\frac{\pi}{2^{k + 1}} = \arctan\left(\frac{\sqrt{2 - c_{k - 1}}}{c_k}\right)
\]
and this completes the proof.
\end{proof}

Since the integer $k$ can be arbitrarily large, we can also write
\begin{equation}
\label{L4AF} 
\frac{\pi }{4} = \lim_{k \to \infty} 2^{k - 1}\arctan\left(\frac{\sqrt{2 - c_{k - 1}}}{c_k}\right).
\end{equation}

Using the limit \eqref{L4AF}, we can derive a well-known formula for $\pi$ \cite{Servi2003} 
\begin{equation}
\label{LWROT} 
\pi = \lim_{k\to\infty} 2^k\sqrt{2 - \underbrace{\sqrt{2 + \sqrt{2 + \sqrt{2 + \cdots +\sqrt{2}}}}}_{k - 1\,\,\rm{square}\,\,\rm{roots}}} = \lim_{k\to\infty} 2^k\sqrt{2 - c_{k - 1}}.
\end{equation}
Another formula for $\pi$ that can also be derived from the limit \eqref{L4AF} is given by (see \cite{Abrarov2025} and literature therein)
\[
\pi = \lim_{k \to \infty}2^k\sum_{n \ge k} \frac{\sqrt{2 - c_{n - 1}}}{c_n}.
\]
It should be noted that this limit can be further simplified as
\[
\pi = \lim_{k \to \infty}2^{k - 1}\sum_{n \ge k} \sqrt{2 - c_{n - 1}}
\]
or
\[
\pi = \lim_{k \to \infty}2^k\sum_{n \ge k} \sqrt{2 - c_n}
\]
since 
\[
\lim_{n \to \infty} c_n = \lim_{n \to \infty} {\underbrace{\sqrt{2 + \sqrt{2 + \sqrt{2 + \cdots +\sqrt{2}}}}}_{n\,\,\rm{square}\,\,\rm{roots}}} = 2.
\]

\section{Irrationality of $\pi$}

\subsection{Motivation}

Define the following constant
\begin{equation}
\label{AK} 
\alpha_k = \left\lfloor{\frac{2^{k + 1}}{\pi}}\right\rfloor,
\end{equation}
where the symbol $\lfloor \cdot \rfloor$ is the floor function that gives the integer part of a number. According to equations \eqref{AF4P} and \eqref{AK}, the constant $\alpha_k$ represents the integer part of the arctangent function as follows:
\[
\alpha_k = \left\lfloor{\frac{1}{\arctan \left(\frac{\sqrt{2 - c_{k - 1}}}{c_k}\right)}}\right\rfloor.
\]
Therefore, we can express the reciprocal of the arctangent function as
\[
\frac{1}{\arctan\left(\frac{\sqrt{2 - c_{k - 1}}}{c_k}\right)} = \alpha_k + \beta_k,
\]
where $\beta_k$ is defined by using the fractional part function $\{ \cdot \}$, as given by
\[
\beta_k = \left\{\frac{1}{\arctan\left(\frac{\sqrt{2 - c_{k - 1}}}{c_k}\right)}\right\} = \left\{\frac{2^{k + 1}}{\pi}\right\}.
\]
Thus, equation \eqref{AF4P} can be expressed in the form
\begin{equation}
\label{AF4PAN} 
\pi  = \frac{2^{k + 1}}{\alpha_k + \beta_k}, \qquad k \ge 1.  
\end{equation}

Since the fractional part $\beta_k$ cannot be smaller than zero and greater than or equal to one while the integer part $\alpha_k$ tends to infinity as $k$ increases, it follows that
\[
\lim_{k \to \infty} \frac{\alpha_k}{\frac{1}{\arctan\left(\frac{\sqrt{2 - c_{k - 1}}}{c_k}\right)}} = \lim_{k \to \infty} \frac{\alpha_k}{\alpha_k + \beta_k} = 1.
\]
Therefore, from this limit and equation \eqref{AF4PAN}, we have
\begin{equation}
\label{L4P} 
\pi  = \lim_{k \to \infty} \frac{2^{k + 1}}{\alpha_k}.
\end{equation}
As we can see from this equation, the integer $\alpha_k$ increases monotonically with increasing $k$.

Consider lemma \ref{L4AK} and theorem \ref{T4NSR} below.
\begin{lemma}
\label{L4AK} 
The following equation
\begin{equation}
\label{IF4A} 
{\alpha _{k + 1}} = \left\{
\begin{aligned}
&2{\alpha _k}, &\qquad 0 \le {\beta_k} < 1/2\\
&2{\alpha _k} + 1, &\qquad 1/2 \le {\beta_k} < 1 & 
\end{aligned}
\right.
\end{equation}
holds.
\end{lemma}

\begin{proof}
From the ratio
\[
\frac{2^{k + 2}\arctan\left( \frac{\sqrt{2 - c_k}}{c_{k + 1}}\right)}{2^{k + 1}\arctan\left(\frac{\sqrt{2 - c_{k - 1}}}{c_k}\right)} = \frac{\pi}{\pi} = 1,
\]
it follows that
\[
\frac{\frac{1}{\arctan \left( \frac{\sqrt{2 - c_k}}{c_{k + 1}}\right)}}{\frac{1}{\arctan\left(\frac{\sqrt{2 - c_{k - 1}}}{c_k}\right)}} = 2
\]
or
\begin{equation}
\label{MR} 
\frac{\alpha_{k + 1} + \beta_{k + 1}}{\alpha_k + \beta_k} = 2.
\end{equation}
Consequently, equation \eqref{MR} results in
\[
\alpha_{k + 1} + \beta_{k + 1} = 2\left(\alpha_k + \beta_k\right)
\]
and with the floor function on both sides, we obtain
\[
\left\lfloor\alpha_{k + 1} + \beta_{k + 1}\right\rfloor  = \left\lfloor 2\left(\alpha_k + \beta_k\right)\right\rfloor.
\]
Since $\alpha_{k + 1}$ is an integer while $0 \le \beta_{k + 1} < 1$, we have
\[
\alpha_{k + 1} = \left\lfloor 2\left(\alpha_k + \beta_k \right) \right\rfloor = 2\alpha_k + \left\lfloor 2\beta_k \right\rfloor
\]
and this completes the proof of lemma \ref{L4AK}.
\end{proof}

\begin{theorem}
\label{T4NSR} 
The following number
\[
\sqrt{2 - {c_k}} , \qquad k \ge 0
\] 
is irrational.
\end{theorem}

\begin{proof}
Suppose that $d$ is an irrational number. It is easy to prove by contradiction that
\begin{equation}
\label{IE} 
\sqrt{2 + d}, \qquad d \in\mathbb{R}\backslash\mathbb{Q}
\end{equation}
is also irrational. In particular, if
\[
\sqrt{2 + d} = \frac{q}{r}, \qquad q,r \in\mathbb{N}\backslash \{0\}
\]
then squaring both sides leads to
\[
2 + d = \left(\frac{q}{r}\right)^2.
\]
Rearranging this equation as
\[
d = \frac{q^2}{r^2} - 2,
\]
results in a contradiction since the right side of this equation is rational. Thus, we have proved that the expression \eqref{IE} must be irrational.

At ${c_0} = 0$, we have
\[
c_1 = \sqrt{2 + c_0}  = \sqrt{2}
\]
and since $\sqrt 2$ is an irrational number \cite{Apostol2000, Ferreno2009} (see also \cite{Dmytryshyn2021}), the next number
\[
{c_2} = \sqrt {2 + {c_1}}  = \sqrt {2 + \sqrt 2 }
\]
is also irrational. The next number 
\[
c_3 = \sqrt{2 + c_2}  = \sqrt{2 + \sqrt{2 + \sqrt{2}}} 
\]
also appears to be irrational since ${c_2}$ is irrational. Repeating this process, by induction we conclude that the number ${c_k}$ must be irrational.

By setting
\[
d = -c_k,
\]
we prove theorem \ref{T4NSR} since the right side of this equation is always irrational.
\end{proof}

Rewrite equation \eqref{LWROT} in the following form:
\[
\pi = \lim_{\ell\to\infty} 2^{\ell + 1}\sqrt{2 - c_\ell}.
\]
We can always choose the integer $\ell$ to be as large as possible such that at any given value of $k$ the approximation
\[
2^{\ell + 1}\sqrt{2 - c_\ell} \approx \pi, \qquad\ell \gg k
\]
retains its accuracy to satisfy the following equation:
\[
\left\lfloor\frac{2^{k + 1}}{2^{\ell + 1}\sqrt{2 - c_\ell}}\right\rfloor = \left\lfloor {\frac{{{2^{k + 1}}}}{\pi }} \right\rfloor, \qquad \;\ell  \gg k.
\]
This equation can be rearranged as
\begin{equation}
\label{EWFF} 
\left\lfloor\frac{2^{k - \ell}}{\sqrt{2 - c_\ell}}\right\rfloor = \left\lfloor \frac{1}{\arctan \left(\frac{\sqrt{2 - c_{k - 1}}}{c_k}\right)}\right\rfloor, \qquad \ell \gg k.
\end{equation}

Consider the following sequence that can be constructed by using either the left or right side of equation \eqref{EWFF}:
\[
\left\{\alpha_k\right\}_{k = 1}^\infty = \{1,2,5,10,20,40,81,162,325,651,1303,2607,5215,10430, \,\,\ldots\}.
\]
The numbers $\alpha_k$ from the sequence $\left\{\alpha_k\right\}_{k = 1}^\infty$ can be found in \cite{OEIS1999a}. Using this sequence, it is convenient to define the following numbers:
\[
\lambda_k = \left\{
\begin{aligned}
0, &\qquad \text{if} \,\, \alpha_k \,\, \text{is even}, \\
1, &\qquad \text{if} \,\, \alpha_k \,\, \text{is odd}.
\end{aligned}
\right.
\]
Therefore, we can construct another sequence:
\begin{equation}
\label{NPS} 
\left\{\lambda_k\right\}_{k = 1}^\infty = \left\{1,0,1,0,0,0,1,0,1,1,1,1,1,0,0,1,1,0,0,0,0,0,1,1,\,\ldots\right\}.
\end{equation}

Under the assumption that the number $\pi$ is rational, we should expect the sequence $\left\{\lambda_k\right\}_{k = 1}^\infty$ to be periodic since any rational number has repeating decimal digits in its expansion. For example, it is well-known that $\pi$ is a number bounded between $3.1408$ and $3.1429$ due to the inequality \cite{Phillips1981, Dalzell1944}:
\[
\frac{223}{71} < \pi < \frac{22}{7}.
\]
Therefore, we can consider the ratio $22/7$ as an example of a rough approximation of $\pi$. In this case, applying
\[
\alpha_k^* = \left\lfloor\frac{2^{k + 1}}{\left(\frac{22}{7}\right)}\right\rfloor 
\]
we can build the following sequence:
\[
\left\{\alpha_k^*\right\}_{k = 1}^\infty = \{1,2,5,10,20,40,81,162,325,651,1303,2606,5213,10426,  \,\ldots\}.
\]
The corresponding counterpart of this sequence, given by
\[
\left\{\lambda_k^*\right\}_{k = 1}^\infty = \left\{1,\underbrace{0,1,0,0,0,1,0,1,1,1}_{\text{10 digits periodicity}},\underbrace{0,1,0,0,0,1,0,1,1,1}_{\text{10 digits periodicity}},\,\ldots\right\},
\]
is not equal to the original sequence $\left\{\lambda_k\right\}_{k = 1}^\infty$ shown by equation \eqref{NPS} above. The periodicity in this sequence occurs because the decimal numbers in this ratio 
\[
\frac{22}{7} = 3.142857142857142857 \ldots  = 3.\overline{142857}
\]
are periodic since it is a rational number.

On the other hand, the sequence $\left\{\lambda_k\right\}_{k = 1}^\infty$ may not be periodic since the expression
\[
\frac{2^{k - \ell}}{\sqrt{2 - c_\ell}}
\]
in equation \eqref{EWFF} is always irrational at any value of $k$ due to the irrationality of the number $\sqrt{2 - c_\ell}$ in accordance with theorem \ref{T4NSR}. This observation with the nested radicals with roots of $2$ strongly motivated us to look for a proof of the irrationality of $\pi$.

It is interesting to note that the sequence \eqref{NPS} coincides with the binary expansion of the constant $1/\pi$ \cite{OEIS1999b}.

\subsection{Bounds and the Dalzell integral}

Define for convenience:
\[
\mu_k = \frac{\sqrt{2 - c_{k - 1}}}{c_k}.
\]

The Maclaurin series expansion of the arctangent function is given by
\begin{equation}
\label{MSE} 
\arctan(x) = \sum_{n = 0}^\infty \frac{(-1)^n x^{2n + 1}}{2n + 1} = x - \frac{x^3}{3} + \frac{x^5}{5} - \frac{x^7}{7} + \cdots \,\,, \qquad |x| \le 1.
\end{equation}
Therefore, in the domain of our interest
\[
x \in \left(0,\mu_k \right],
\]
where
\[
0 < \mu_k \le 1, \qquad k \ge 1,
\]
the difference between $x$ and $\arctan(x)$ can be squeezed between the boundaries $0$ and ${x^3}/3$ as shown by the following inequality:
\begin{equation}
\label{BI} 
0 < x - \arctan(x) < \frac{x^3}{3}.
\end{equation}

Substituting $x = \mu_k$, we can rewrite inequality \eqref{BI} as
\[
0 < \mu_k - \arctan(\mu_k) < \frac{\mu_k^3}{3}
\]
or
\begin{equation}
\label{IWMK} 
0 < 2^{k + 1}\mu_k - 2^{k + 1}\arctan(\mu_k) < \frac{2^{k + 1}\mu_k^3}{3}.
\end{equation}

According to the equation
\[
\mu_k - \arctan(\mu_k) = \int_0^{\mu_k} \frac{x^2}{1 + x^2}dx
\]
the difference ${\mu _k} - \arctan(\mu_k)$ can be interpreted geometrically as the area under the curve $x^2/\left(1 + x^2\right)$ bounded between points $0$ and $\mu_k$. More generally, a relation with the Maclaurin series expansion \eqref{MSE} can also be shown as
\[
\sum_{n = 0}^m \frac{(-1)^n \mu_k^{2n + 1}}{2n + 1} - \arctan(\mu_k) = -\int_0^{\mu_k} \frac{(-1)^{m + 1} x^{2(m + 1)}}{1 + x^2}dx.
\]

Taking into consideration that
\[
2^{k + 1}\arctan(\mu_k) = \pi
\]
we can express inequality \eqref{IWMK} in the form
\begin{equation}
\label{IWP} 
0 < 2^{k + 1} \mu_k - \pi  < \frac{2^{k + 1} \mu_k^3}{3}.
\end{equation}

Assume that the constant $\pi$ is rational and can be represented as a ratio of two positive integers $q$ and $r$. Under this assumption, inequality \eqref{IWP} can be rewritten as
\[
0 < 2^{k + 1} \mu_k - \frac{q}{r} < \frac{2^{k + 1}\mu _k^3}{3}
\]
or
\begin{equation}
\label{IWI} 
0 < 2^{k + 1}\mu_k r - q < \frac{2^{k + 1} \mu_k^3 r}{3}.
\end{equation}

Since $r$ is a positive integer, the number $2^{k + 1}\mu_k r$ is also positive. Therefore, this number can be expanded in terms of its integer and fractional parts as
\[
2^{k + 1}\mu_k r = \left\lfloor 2^{k + 1}\mu_k r \right\rfloor + \left\{2^{k + 1}\mu_k r\right\},
\]
where
\begin{equation}
\label{I4M} 
0 \le \{2^{k + 1}\mu_k r\} < 1, \qquad k \ge 1.
\end{equation}

Substituting the last equation above into inequality \eqref{IWI}, we get
\[
0 < \left\lfloor 2^{k + 1}\mu_k r \right\rfloor + \left\{2^{k + 1}\mu_k r\right\} - q < \frac{2^{k + 1}\mu_k^3 r}{3}
\]
or
\begin{equation}
\label{FI} 
-\left\{2^{k + 1}\mu_k r\right\} < \left\lfloor 2^{k + 1}\mu_k r \right\rfloor - q < \frac{2^{k + 1}\mu_k^3 r}{3} - \left\{2^{k + 1}\mu_k r\right\}.
\end{equation}

Since
\[
\frac{1}{\mu_k} = \frac{c_k}{\sqrt{2 - c_{k - 1}}} \ge \alpha_k \ge 2^{k - 1}
\]
we can determine that
\[
\frac{2^{k + 1}\mu_k^3 r}{3} = \frac{4r}{3}2^{k - 1}\mu_k^3 = \frac{4r}{3}\frac{2^{k - 1}}{1/\mu_k^3} \le \frac{4r}{3}\mu_k^2
\]
and since
\[
\frac{4r}{3}\mu_k^2 = \frac{4r}{3}\left(\frac{\sqrt{2 - c_{k - 1}}}{c_k}\right)^2 \to 0, \qquad k \to \infty
\]
we end up with
\[
\frac{2^{k + 1}\mu_k^3 r}{3} \to 0, \qquad k \to \infty.
\]

It may be interesting to demonstrate empirical evidence supporting the irrationality of $\pi$ by attempting to lock the integer $\left\lfloor 2^{k + 1}r\mu_k\right\rfloor - q$ between the lower and upper bounds of inequality \eqref{FI} at a sufficiently large $k$. In particular, we can observe that at given $q$ and $r$ we can find some large integer $k$ for which (and above which) inequality \eqref{FI} will no longer be valid. This occurs because at a sufficiently large $k$ there will not be enough space between the lower
\[
-\left\{2^{k + 1}\mu_k r\right\}
\]
and upper
\[
\frac{2^{k + 1}\mu_k^3 r}{3} - \left\{2^{k + 1}\mu_k r\right\}
\]
bounds to accommodate the integer
\[
\left\lfloor 2^{k + 1}\mu_k r\right\rfloor - q
\]
inside the domain
\[
D = \left(-\left\{2^{k + 1}\mu_k r\right\},\frac{2^{k + 1}\mu_k^3 r}{3} - \left\{2^{k + 1}\mu_k r\right\}\right)
\]
that very rapidly narrows with an increasing $k$.

Consider a few examples by generating rational approximations of $\pi$ based on the (generalized) Dalzell integral \cite{Nield1982, Lucas2009} (see also \cite{Dalzell1944, Backhouse1995, Shiu2014, Rattaggi2018}):
\[
\begin{aligned}
&\int_0^1 \frac{x^{4m}(1 - x)^{4m}}{1 + x^2}dx = \\
&-(-4)^{m - 1}\pi + \int_0^1 \left(x^6 - 4x^5 + 5x^4 - 4x^2 + 4\right)\sum_{n = 0}^{m - 1} (-4)^{m - 1 - n}x^{4n}(1 - x)^{4n}dx.  
\end{aligned}
\]
Since the Dalzell integral above vanishes with an increasing integer $m$, we can use it to approximate $\pi$ as follows:
\begin{equation}
\label{A4P} 
\begin{aligned}
\pi \approx &\frac{1}{(-4)^{m - 1}} \\
&\times\int_0^1 \left(x^6 - 4x^5 + 5x^4 - 4x^2 + 4\right)\sum_{n = 0}^{m - 1}(-4)^{m - 1 - n}x^{4n}(1 - x)^{4n}dx.
\end{aligned}
\end{equation}

At $m = 1$, the approximation \eqref{A4P} yields \cite{Dalzell1944}
\[
\pi \approx \frac{22}{7}.
\]
Substituting $q = 22$ and $r = 7$ into inequality \eqref{FI}, we can find that it remains valid only while $k \le 5$. At $k = 6$, inequality \eqref{FI} becomes invalid since the corresponding integer
\[
\left\lfloor 2^{6 + 1}\mu_6 r \right\rfloor - q = -1,
\]
is not within the lower
\[
-\left\{2^{6 + 1}\mu_6 r\right\} = -0.9955654095 \ldots
\]
and upper
\[
\frac{2^{6 + 1}\mu_6^3 r}{3} - \left\{2^{6 + 1}\mu_6 r\right\} = -0.9911469781 \ldots 
\]
bounds. In particular, at $m = 1$ and $k = 6$ it appears that
\[
-\left\{2^{6 + 1}\mu_6 r\right\} > \left\lfloor 2^{6 + 1}\mu_6 r\right\rfloor - q
\]
and this contradicts inequality \eqref{FI}.

At $m = 4$, the approximation \eqref{A4P} provides
\[
\pi \approx \frac{741,269,838,109}{235,953,517,800}.
\]
Substituting $q = 741,269,838,109$ and $r = 235,953,517,800$ into inequality \eqref{FI}, we can see that it remains valid only until $k \le 10$. At $k = 11$, inequality \eqref{FI} is no longer valid since the corresponding integer
\[
\left\lfloor 2^{11 + 1}\mu_{11}r\right\rfloor - q = 145,356
\]
appears to be beyond the lower
\[
-\left\{2^{11 + 1}\mu_{11}r\right\} = -0.8389619580 \ldots
\]
and upper
\[
\frac{2^{11 + 1}\mu_{11}^3 r}{3} - \left\{2^{11 + 1}\mu_{11}r\right\} = 145,355.9027929868 \ldots
\]
bounds. Specifically, at $m = 4$ and $k = 11$ we can find that
\[
\left\lfloor 2^{11 + 1}\mu_{11}r\right\rfloor - q > \frac{2^{11 + 1}\mu_{11}^3 r}{3} - \left\{2^{11 + 1}\mu_{11}r\right\}
\]
contradicts inequality \eqref{FI}.

At $m = 8$, the approximation \eqref{A4P} provides
\[
\pi \approx \frac{19,809,071,774,292,917,047,896,724,979}{6,305,423,381,881,718,760,060,595,200}.
\]
Substituting
\[
q = 19,809,071,774,292,917,047,896,724,979
\]
and
\[
r = 6,305,423,381,881,718,760,060,595,200
\]
into inequality \eqref{FI}, we can determine that it remains valid only for $k \le 20$. At $k = 21$, inequality \eqref{FI} becomes invalid since the corresponding integer
\[
\left\lfloor 2^{21 + 1}\mu_{21}r\right\rfloor - q = 3,704,442,064,441,320
\]
is not inside the lower
\[
-\left\{2^{21 + 1}\mu_{21}r\right\} = -0.0467121698 \ldots
\]
and upper
\[
\frac{2^{21 + 1}\mu_{21}^3 r}{3} - \left\{2^{21 + 1}\mu_{21}r\right\} = 3,704,442,064,439,983.9322317187 \ldots
\]
bounds. In particular, at $m = 8$ and $k = 21$ we can see that
\[
\left\lfloor 2^{21 + 1}\mu_{21}r\right\rfloor - q > \frac{2^{21 + 1}\mu_{21}^3 r}{3} - \left\{2^{21 + 1}\mu_{21}r\right\}
\]
contradicts inequality \eqref{FI} again.

From these examples, we can observe the tendency indicating that increasing the accuracy of $\pi$ by using a rational approximation $q/r$ with larger and larger integers in its numerator $q$ and denominator $r$ increases the range of the integer $k$ at which inequality \eqref{FI} remains valid. Therefore, this observation is consistent with the irrationality of $\pi$.

Consider now $m = 8$ when $k \gg 1$. We can take, say, $k = 100$. Since in this case
\[
\frac{2^{100 + 1}\mu_{100}^3 r}{3} = 1.01387 \ldots \times 10^{-32} \ll 1,
\]
we can observe that the left
\[
-\left\{2^{100 + 1}\mu_{100}r\right\} = -0.03199241354446692964676896442009338 \ldots 
\]
and right
\[
\begin{aligned}
\frac{2^{100 + 1}\mu_{100}^3 r}{3} &- \left\{2^{100 + 1}\mu_{100}r\right\} \\
&= -0.03199241354446692964676896442008324 \ldots
\end{aligned}
\]
bounds are nearly the same. Therefore, when $k \gg 1$, due to a very narrow gap between the left and right bounds in inequality \eqref{FI}, even a single integer cannot exist inside them. This empirical evidence is also consistent with the irrationality of $\pi$.

\subsection{Continued fractions}

Consider the Śleszyński--Pringsheim theorem \ref{SPT} below \cite{Lorentzen2018}.

\begin{theorem}
\label{SPT}
Let $a_n, b_n \in \mathbb{R}$ for all $n \ge 1$ with $a_n \neq 0$. If
$$
\vert{}b_n\vert{} \ge \vert{}a_n\vert{} + 1 \quad \text{for all} \,\,n \ge 1,
$$
then the sequence of approximants $g_n$ of the continued fraction
$$
g = \cfrac{a_1}{b_1 + \cfrac{a_2}{b_2 + \cfrac{a_3}{b_3 + \dots}}}
$$
converges to a real limit
$$
\lim_{n\to\infty} g_n = g \in \mathbb{R}.
$$
Moreover, every approximant satisfies
$$
\vert{}g_n\vert{} < 1,
$$
and the limit satisfies 
$$
\lim_{n\to\infty} \vert{}g_n\vert{} = \vert{}g\vert{} \le 1.
$$
\end{theorem}

\begin{proof}
The $n$-th approximant is given by $g_n = A_n/B_n$ for $n \ge 1$, where $A_n$ and $B_n$ can be determined by the following recurrence relations (two-step iteration):
$$
A_n = b_n A_{n-1} + a_n A_{n-2} \\
$$
\begin{equation}
\label{IF4BN} 
B_n = b_n B_{n-1} + a_n B_{n-2},
\end{equation}
where initial values are $A_{-1} = 1, A_0 = 0, B_{-1} = 0$ and $B_0 = 1$.

Applying the reverse triangle inequality to equation \eqref{IF4BN}, we obtain
$$
\vert{}B_n\vert{} \ge \vert{}b_n\vert{} \vert{}B_{n-1}\vert{} - \vert{}a_n\vert{} \vert{}B_{n-2}\vert{}.
$$
Using now the condition $\vert{}b_n\vert{} \ge \vert{}a_n\vert{} + 1$, we get
$$
\vert{}B_n\vert{} \ge (\vert{}a_n\vert{} + 1) \vert{}B_{n-1}\vert{} - \vert{}a_n\vert{} \vert{}B_{n - 2}\vert{}
$$
or
$$
\vert{}B_n\vert{} \ge \vert{}B_{n - 1}\vert{} + \vert{}a_n\vert{} (\vert{}B_{n - 1}\vert{} - \vert{}B_{n-2}\vert{}).
$$
Subtracting $\vert{}B_{n - 1}\vert{}$ from both sides gives the following step-wise inequality:
$$
\vert{}B_n\vert{} - \vert{}B_{n - 1}\vert{} \ge \vert{}a_n\vert{} (\vert{}B_{n - 1}\vert{} - \vert{}B_{n - 2}\vert{}).
$$

Thus, for the case $n = 1$, we have
$$
\vert{}B_1\vert{} - \vert{}B_0\vert{} = \vert{}b_1\vert{} - 1 \ge \vert{}a_1\vert{}.
$$
Consequently, by induction we obtain
$$
\vert{}B_n\vert{} - \vert{}B_{n-1}\vert{} \ge \prod_{m = 1}^n \vert{}a_m\vert{} > 0.
$$

The inequality $\vert{}a_m\vert{} > 0$ implies that $\{\vert{}B_n\vert{}\}_{n = 0}^\infty$ is a strictly increasing sequence, where $\vert{}B_0\vert{} = 1$. Therefore, we can infer that
$$
\vert{}B_n\vert{} > 1, \qquad n \ge 1.
$$

Using the determinant formula for continued fractions, we obtain
\begin{equation}
\label{PF} 
A_n B_{n - 1} - A_{n - 1} B_n = (-1)^{n - 1} \prod_{m = 1}^n a_m.
\end{equation}
Consequently, setting $g_0 = 0$, we can find the difference between consecutive values of the approximants $g_n - g_{n-1}$. Thus, dividing both sides of equation \eqref{PF} by $B_n B_{n-1}$ yields
$$
g_n - g_{n - 1} = \frac{A_n}{B_n} - \frac{A_{n - 1}}{B_{n - 1}} = \frac{(-1)^{n - 1} \prod_{m = 1}^n a_m}{B_n B_{n - 1}}.
$$
Therefore, the absolute value for the difference is given by
$$
\vert{}g_n - g_{n - 1}\vert{} = \frac{\prod_{m = 1}^n \vert{}a_m\vert{}}{\vert{}B_n B_{n - 1}\vert{}}.
$$

Substituting our lower bound 
$$
\prod_{m = 1}^n \vert{}a_m\vert{} \le \vert{}B_n\vert{} - \vert{}B_{n - 1}\vert{}
$$
into the numerator
$$
\vert{}g_n - g_{n - 1}\vert{} \le \frac{\vert{}B_n\vert{} - \vert{}B_{n - 1}\vert{}}{\vert{}B_n\vert{} \vert{}B_{n - 1}\vert{}} = \frac{1}{\vert{}B_{n - 1}\vert{}} - \frac{1}{\vert{}B_n\vert{}}
$$
and summing these differences from $m = 1$ to $n$, results in
$$
\sum_{m = 1}^n \vert{}g_m - g_{m - 1}\vert{} \le \sum_{m = 1}^n \left( \frac{1}{\vert{}B_{m - 1}\vert{}} - \frac{1}{\vert{}B_m\vert{}} \right) = \frac{1}{\vert{}B_0\vert{}} - \frac{1}{\vert{}B_n\vert{}} = 1 - \frac{1}{\vert{}B_n\vert{}} < 1.
$$

The partial sums of the form $\sum_{m = 1}^n \vert{}g_m - g_{m - 1}\vert{}$ are bounded above by $1$. Therefore, we can conclude that the series $\sum_{m = 1}^\infty (g_m - g_{m - 1})$ is absolutely convergent. Since absolute convergence implies the convergence of the series, the sequence of approximants $\{g_n\}_{n = 0}^\infty$ can be determined by using the following partial sums:
$$
g_n = \sum_{m = 1}^n (g_m - g_{m - 1}).
$$
Therefore, $g_n$ converges to a real limit such that
$$
g = \lim_{n \to \infty} g_n \in \mathbb{R}.
$$
Furthermore, due to the triangle inequality we can write
$$
\vert{}g_n\vert{} \le \sum_{m = 1}^n \vert{}g_m - g_{m - 1}\vert{} < 1
$$
from which it follows that
$$
\vert{}g_n\vert{} < 1 \quad \text{and} \quad \vert{}g\vert{} \le 1
$$
and this completes the proof of the theorem.
\end{proof}

The next important theorem \ref{T4CIT} is also related to infinite continued fractions.

\begin{theorem}
\label{T4CIT} 
Let $d\in 2\mathbb Z\setminus\{0\}$. Define
$$
x_N = 1 + \underset{n = 1}{\overset{N}{\mathbf K}}\frac{d}{2n + 1} = 1 + \cfrac{d}{3 + \cfrac{d}{5 + \cfrac{d}{\ddots + \cfrac{d}{2N + 1}}}}, \qquad N \ge 1,
$$
where $\mathbf{K}$ denotes the Kettenbruch notation (also known as Gauss’s K-notation). Then, the sequence $\{x_N\}_{N = 1}^{\infty}$ converges to a finite real limit
$$
x_\infty = \lim_{N\to\infty} x_N \in \mathbb{R}.
$$
\end{theorem}

\begin{proof}
Let
$$
x_N = \frac{P_N}{Q_N},
$$
where the sequences $\{P_n\}_{n = -1}^{\infty}$ and $\{Q_n\}_{n = -1}^{\infty}$ are defined by
$$
\begin{aligned}
P_n &= (2n + 1)P_{n - 1} + dP_{n - 2},\\
Q_n &= (2n + 1)Q_{n - 1} + dQ_{n - 2},
\end{aligned}
\qquad n \ge 1,
$$
with initial values
$$
P_{-1} = 1, \qquad P_0 = 1,
$$
and
$$
Q_{-1} = 0, \qquad Q_0 = 1.
$$
These equations are the standard continuant recurrences (or continuant iterations) associated with the finite continued fraction.

We first observe that every $Q_n$ is odd. Indeed, since $Q_0 = 1$ is odd, mathematical induction yields
$$
Q_n = (2n + 1)Q_{n - 1} + d Q_{n - 2} \equiv 1 \pmod2,
$$
since $2n + 1$ is odd and $d$ is even. Therefore, this gives
$$
Q_n \equiv 1 \pmod2, \qquad n \ge 0.
$$
Consequently, since $Q_n \ne 0$ for all $n \ge 0$, every approximant $x_N$ is finite.

We can choose
$$
M_1 \ge \frac{|d|}{2}
$$
and define the finite tail as
$$
t_{M_1,N} = \underset{j = 0}{\overset{N - M_1}{\mathbf K}} \frac{d}{2(M_1 + j) + 1}, \qquad N \ge M_1.
$$
The corresponding partial numerators and denominators of this tail are
$$
a_j = d, \quad b_j = 2(M_1 + j) + 1.
$$

Since $2M_1 \ge |d|$, we have
$$
|b_j| = 2(M_1 + j) + 1 \ge 2M_1 + 1 \ge |d| + 1 = |a_j| + 1.
$$

Consequently, by the Śleszyński--Pringsheim theorem (\ref{SPT}),
$$
t_{M_1,N} \rightarrow L
$$
for some finite real number $L$ satisfying $|L| \le 1$.

We next establish the eventual strict growth of the absolute values of the continuant denominators. Since every $Q_n$ is a nonzero integer, we have
$$
|Q_n| \ge 1, \qquad n \ge 0.
$$

We claim there exists an arbitrarily large index $m$ such that
$$
|Q_m| > |Q_{m - 1}|.
$$
Suppose for contradiction this is false. Then, there would exist an index $N_0$ such that
$$
|Q_n| \le |Q_{n-1}|, \qquad n \ge N_0.
$$
Hence, $\{|Q_n|\}_{n\ge N_0 - 1}$ would be a non-increasing sequence of positive integers, meaning it must eventually become constant. Thus, for some integer $B \ge 1$, we get
$|Q_n| = B$ for all sufficiently large $n$. However, applying the triangle inequality to the recurrence relation yields
$$
(2n + 1)|Q_{n - 1}| \le |Q_n| + |d|\,|Q_{n - 2}|.
$$
Since $|Q_{n}| = |Q_{n-1}| = |Q_{n-2}| = B$ for all sufficiently large $n$, this implies that
$$
(2n + 1)B \le (1 + |d|)B,
$$
which is impossible for sufficiently large $n$. This contradiction proves our claim.

We can choose such an index $m$ that is large enough to satisfy
$$
2m + 1 - |d| \ge 3,
$$
and set $M_0 = m + 1$. Then, by our choice of $m$, we can write
$$
|Q_{M_0 - 1}| > |Q_{M_0 - 2}| > 0.
$$

We are now ready to prove by induction that
$$
|Q_n| > |Q_{n-1}| > 0, \qquad n \ge M_0.
$$

Assume that
$$
|Q_{n - 1}| > |Q_{n - 2}| > 0
$$
for some $n \ge M_0$. From the recurrence relation and the reverse triangle inequality, we obtain
$$
\begin{aligned}
|Q_n| &\ge (2n + 1)|Q_{n - 1}| - |d|\,|Q_{n - 2}| \\
&> (2n + 1)|Q_{n - 1}| - |d|\,|Q_{n - 1}| \\
&= (2n + 1 - |d|)|Q_{n - 1}|.
\end{aligned}
$$

Since $n \ge M_0$, it follows that $2n + 1 - |d| \ge 3$. Therefore, the inequality becomes
$$
|Q_n| > 3|Q_{n - 1}| > |Q_{n-1}| > 0,
$$
which completes the induction step.

Define
$$
M = \max\left(M_0, M_1\right).
$$
Then, for every $N \ge M$, the continuant decomposition gives
$$
x_N = \frac{P_{M - 1} + t_{M,N}P_{M - 2}}{Q_{M - 1} + t_{M,N}Q_{M - 2}},
$$
where
$$
t_{M,N} = \underset{k=0}{\overset{N-M}{\mathbf K}}\frac{d}{2(M + k) + 1}.
$$

Since $M \ge M_1$, the Śleszyński--Pringsheim theorem \ref{SPT} guarantees that
$$
t_{M,N} \rightarrow L_M
$$
for some $L_M \in \mathbb{R}$ satisfying $|L_M| \le 1$.

Define the limiting denominator as
$$
D_M = Q_{M-1} + L_M Q_{M-2}.
$$

Because $M \ge M_0$, we know that $|Q_{M-1}| > |Q_{M-2}| \ge 0$. Using the reverse triangle inequality and the fact that $|L_M| \le 1$, we find
$$
\begin{aligned}
|D_M| &\ge |Q_{M-1}| - |L_M|\,|Q_{M-2}| \\
&\ge |Q_{M-1}| - |Q_{M-2}| \\
&> 0.
\end{aligned}
$$
Thus, we show that $D_M \ne 0.$

Since $t_{M,N} \rightarrow L_M$ and the denominator of the limiting Möbius transformation \cite{Schwerdtfeger1979} is strictly non-zero, we can write
$$
\begin{aligned}
x_N &= \frac{P_{M-1} + t_{M,N}P_{M-2}}{Q_{M-1} + t_{M,N}Q_{M-2}} \\
&\longrightarrow \frac{P_{M-1} + L_M P_{M-2}}{Q_{M-1} + L_M Q_{M-2}}.
\end{aligned}
$$
Therefore, we can conclude that
$$
\boxed{
x_\infty = \frac{P_{M-1} + L_M P_{M-2}}{Q_{M-1} + L_M Q_{M-2}} \in \mathbb{R}
}
$$
and this completes the proof that the sequence $\{x_N\}_{N=1}^{\infty}$ converges to a finite real limit.
\end{proof}

We can show now how the irrationality of the constant $\pi$ can be proved by using theorem \ref{T4CIT} above. In particular, under the assumption that $\pi$ is rational, it is straightforward to observe that
$$
\tan(r\pi) = 0 \implies \tan\left(r\frac{q}{r}\right) = \tan(q) = 0.
$$
Therefore, it is sufficient to prove that
\[
\tan(q) \ne 0
\]
to prove the irrationality of $\pi$. In this context, the statement "the tangent of a positive integer cannot be zero" can be regarded as a reformulation of the irrationality of $\pi$. This is shown in the next theorem \ref{T4IP}.

\begin{theorem}
\label{T4IP} 
The number $\pi$ is irrational.
\end{theorem}

\begin{proof}
Suppose that the constant $\pi$ can be represented as a ratio of two integers $q$ and $r$ such that
$$
\pi = \frac{q}{r}, \qquad q,r \in \mathbb{N}\backslash\{0\}.
$$
Then, from equation \eqref{CF4TF} it follows that
\begin{equation}
\label{TF4I} 
\tan(r\pi) = \tan(q) = \frac{q}{y_\infty},
\end{equation}
where
$$
y_\infty = 1 + \cfrac{-q^2}{3 + \cfrac{-q^2}{5 + \cfrac{-q^2}{7 + \cfrac{-q^2}{9 + \ddots}}}}.
$$
We will consider two cases separately for an even and an odd integer $q$.

If we assume that $q$ is an even integer, then from theorem \ref{T4CIT} it follows that
$$
\vert{}y_\infty\vert{} < \infty.
$$
As a result, equation \eqref{TF4I} cannot be equal to zero since its denominator $y_{\infty}$ is convergent.

Assume now that $q$ is an odd integer. In this case, we cannot use theorem \ref{T4CIT} directly. However, using the following equation
$$
\tan(2r\pi) = \tan(2q)
$$
and considering the fact that $2q$ is an even integer, we can apply theorem \ref{T4CIT} again. Thus, representing the tangent function by a continued fraction as
\begin{equation}
\label{TF4DQ} 
\tan(2q) = \frac{2q}{z_\infty},
\end{equation}
where
$$
z_\infty = 1 + \cfrac{-(2q)^2}{3 + \cfrac{-(2q)^2}{5 + \cfrac{-(2q)^2}{7 + \cfrac{-(2q)^2}{9 + \ddots}}}},
$$
we come to the conclusion that according to theorem \ref{T4CIT}, equation \eqref{TF4DQ} cannot be equal to zero since
$$
\vert{}z_\infty\vert{} < \infty.
$$

Under the assumption that $\pi$ is rational, from the identity
$$
\tan(r\pi) = \tan(2r\pi) = 0,
$$
it follows that
$$
\tan(q) = \tan(2q) = 0.
$$
However, we have shown already that $\tan(2q) \ne 0$. Therefore, the value $\tan(q)$ also cannot be equal to zero when $q$ is odd.

Since the tangent function of any even or odd integer $q$ cannot be equal to zero, we have a contradiction with the equation $\tan(r\pi) = \tan(q) = 0$. Therefore, the constant $\pi$ cannot be represented as a ratio of two integers $q$ and $r$. This completes the proof that the constant $\pi$ is irrational.
\end{proof}

Theorem \ref{T4IP} can also be proved by using the nested radicals with square roots of 2. Specifically, since from theorem \ref{T4CIT} it follows that $\tan(2q)$ cannot be equal to zero, the following equation
\[
\lim_{k\to\infty} \tan\left(2^{k + 2}r\tan\left(\frac{\pi}{2^{k + 1}}\right)\right) = \lim_{k\to\infty} \tan\left(2^{k + 2}r\mu_k\right) = \tan(2r\pi) = 0,
\]
leads to
\[
\lim_{k\to\infty} 2^{k + 2}r\mu_k \ne 2q
\]
or
\[
\lim_{k\to\infty} 2^{k + 1}r\mu_k \ne q,
\]
from which it follows that
\[
\left\{\lim_{k\to\infty} 2^{k + 1}r\mu_k\right\} \ne 0 \implies \left\{\lim_{k\to\infty} 2^{k + 1}r\mu_k\right\} = \{r\pi\} > 0
\]
proving the irrationality of $\pi$.

It should be noted that Lambert in his proof used a strictly decreasing sequence of positive integers ultimately leading to a contradiction with the irrationality of $\pi$ \cite{Berggren2004, Angell2022}. However, in this approach, we proved by using the Śleszyński--Pringsheim theorem \ref{SPT} that the tangent of any non-zero even integer cannot be equal to zero. Since this methodology is conceptually different from that used by Lambert, the proof of the irrationality of $\pi$ presented here is not a modification of original Lambert's proof.

\subsection{Irrationality of $\mu_{k \ge 2}$}

We can prove that the number $\tan\left(\pi/2^{k + 1}\right) = \mu_k$ is also irrational at $k > 1$. The proof is straightforward. Consider theorem \eqref{T4M} below.
\begin{theorem}
\label{T4M} 
At $k \ge 2$ the number $\mu_k$ is irrational.
\end{theorem}

\begin{proof}
Suppose that the number $\tan(\pi/16)$ is rational and, therefore, can be written such that
\[
\tan\left(\frac{\pi}{16}\right) = \tan\left(\frac{\pi}{2^{3 + 1}}\right) = \mu_3 = \frac{\sqrt{2 - \sqrt{2 + \sqrt{2}}}}{\sqrt{2 + \sqrt{2 + \sqrt{2}}}} = \frac{q}{r}, \qquad q,r \in\mathbb{N}\backslash\{0\}.
\]
Then, according to the double-angle formula for the tangent function, the number
\[
\tan\left(\frac{\pi}{8}\right) = \tan\left(\frac{\pi}{2^{2 + 1}}\right) = \frac{2\tan(\pi/16)}{1 - \tan^2(\pi/16)}
\]
is also rational since the numerator $2\tan(\pi/16)$ and denominator $1 - \tan^2(\pi/16)$ are both rational numbers. However, we know that
\[
\tan\left(\frac{\pi}{8}\right) = \mu_2 = \frac{\sqrt{2 - \sqrt{2}}}{\sqrt{2 + \sqrt{2}}} = \sqrt{2} - 1
\]
is irrational since $\sqrt{2}$ is irrational \cite{Apostol2000, Ferreno2009} and we reached a contradiction. Thus, because of this contradiction, we prove that the number $\tan(\pi/16)$ is irrational. Since $\tan(\pi/16)$ is an irrational number, we can prove now that
\[
\tan\left(\frac{\pi}{32}\right) = \tan\left(\frac{\pi}{2^{4 + 1}}\right) = \mu_4 = \frac{\sqrt{2 - c_{4 - 1}}}{c_4} = \frac{\sqrt{2 - \sqrt{2 + \sqrt{2 + \sqrt{2}}}}}{\sqrt{2 + \sqrt{2 + \sqrt{2 + \sqrt{2}} }}}
\]
is also irrational by using the double-angle formula for the tangent function again. Thus, by telescoping this induction procedure further for larger and larger values of $k$, we conclude the proof that
\[
\mu_k \in \mathbb{R}\setminus\mathbb{Q}, \qquad k \ge 2.
\]
\end{proof}

As $\mu_{k \ge 2}$ is irrational, we can modify inequality \eqref{I4M} as
\[
0 < \{2^{k + 1}\mu_k r\} < 1, \qquad k \ge 2.
\]


\section{Rational approximation of $\pi$}

The limit \eqref{L4P} can be used to generate a rational approximation of $\pi$ and we can show the importance of the odd numbers from the sequence $\left\{\alpha_k\right\}_{k = 1}^\infty$ to clarify how the rational approximation approaches the constant $\pi$ as $k$ increases.

Define the following integers:
\[
\gamma_k = \left\{
\begin{aligned}
&k, & & {\text{if}} \,\, \alpha_k \,\, \text{is odd}, \\
&\gamma_{k - 1}, & & {\text{if}} \,\, \alpha_k \,\, \text{is even}.
\end{aligned}
\right.
\]
Now we can construct the sequences for positive integers $\gamma_k$ and $\alpha_{\gamma_{_k}}$ as follows:
\[
\left\{\gamma_k\right\}_{k = 1}^\infty = \{1,1,3,3,3,3,7,7,9,10,\,\ldots\}
\]
and
\[
\begin{aligned}
\left\{\alpha_{\gamma_{_k}}\right\}_{k = 1}^\infty &= \left\{\alpha_1,\alpha_1,\alpha_3,\alpha_3,\alpha_3,\alpha_3,\alpha_7,\alpha_7,\alpha_9,\alpha_{10},\,\,\ldots\right\} \\
&= \{1,1,5,5,5,5,81,81,325,651,\,\ldots\}.
\end{aligned}
\]

According to lemma \ref{L4AK}, the integers in the sequence $\left\{\alpha_k\right\}_{k = 1}^\infty$ can be even and odd. However, the integers in the sequence $\left\{\alpha_{\gamma_{_k}}\right\}_{k = 1}^\infty$ are always odd. This means that if an integer $\alpha_k$ is an even number, then it has a common factor with the integer $2^{k + 1}$ as both of them are divisible by $2$. Thus, due to divisibility by $2$ when $\alpha_k$ is an even number, we can rearrange the limit \eqref{L4P} as
\begin{equation}
\label{ML} 
\pi = \lim_{k \to \infty} \frac{2^{\gamma_{_k} + 1}}{\alpha_{\gamma_{_k}}}, \qquad \alpha_{\gamma_{_k}} \in 2\mathbb{N} + 1.
\end{equation}

The limits \eqref{L4P} and \eqref{ML} show that we can approximate $\pi$ in the form of a rational approximation as given by
\begin{equation}
\label{RA4P} 
\pi \approx \frac{2^{k + 1}}{\alpha_k} = \frac{2^{\gamma_k + 1}}{\alpha_{\gamma_{_k}}}, \qquad k \gg 1.
\end{equation}
This means that numerator $2^{k + 1}$ and denominator $\alpha_k$ cannot be reduced smaller than $2^{\gamma_k + 1}$ and $\alpha_{\gamma_{_k}}$, respectively, since even $2^{\gamma_k + 1}$ and odd $\alpha_{\gamma_{_k}}$ are always relatively prime numbers that have no common divisors other than $1$.

Definition \eqref{IF4A} implies that the fractional part
\[
0 \le \beta_k < 1.
\]
However, we can also show that the following strict inequality
\[
0 < \beta_k < 1
\]
is also valid.

The fractional $\beta_k$ is always greater than zero because the reciprocal of the arctangent function in the following equation
\[
\frac{1}{{\arctan \left( {\frac{{\sqrt {2 - {c_{k - 1}}} }}{{{c_k}}}} \right)}} = \frac{{{2^{k + 1}}}}{\pi}
\]
cannot be an integer due to the irrationality of $\pi$. Consequently, the fractional part
(see equation \eqref{AF4PAN})
\[
\beta_k = \frac{{{2^{k + 1}}}}{\pi} - \alpha_k, \qquad k \ge 1
\]
cannot be equal to zero.

It is not difficult to show that due to the condition $\beta_k > 0$ it follows that the rational approximation
\[
\frac{2^{\gamma_k + 1}}{\alpha_{\gamma_{_{k}}}} > \pi.
\]
More specifically, the left part of this inequality represents an upper bound of the constant $\pi$. It is also not difficult to show that due to the condition $\beta_k < 1$ the left side of the following inequality
\[
\frac{2^{\gamma_k + 1}}{\alpha_{\gamma_{_{k}}} + 1} < \pi
\]
provides a lower bound of the constant $\pi$. For instance, by taking $k = 10$ we obtain
\[
\frac{512}{163} < \pi < \frac{2048}{651}.
\]

Consider now the following examples (a link for the extended table showing values of $\alpha_k$ can be found in \cite{OEIS1999a}):
\[
\begin{aligned}
&\alpha_{70} & & = \alpha_{\gamma_{_{70}}} & & = 751,587,968,840,192,313,983 \\
&\alpha_{71} & & = 2\alpha_{\gamma_{_{70}}} & & = 1,503,175,937,680,384,627,966 \\
&\alpha_{72} & & = 4\alpha_{\gamma_{_{70}}} & & = 3,006,351,875,360,769,255,932 \\
&\alpha_{73} & & = 8\alpha_{\gamma_{_{70}}} & & = 6,012,703,750,721,538,511,864 \\
&\alpha_{74} & & = 16\alpha_{\gamma_{_{70}}} & & = 12,025,407,501,443,077,023,728
\end{aligned}
\]
Although the values of the coefficients from ${\alpha _{70}}$ to ${\alpha _{74}}$ increase by a factor of $2$, the corresponding ratios
\[
\begin{aligned}
\frac{2^{75}}{\alpha_{74}} &= \frac{2^{74}}{\alpha _{73}} = \frac{2^{73}}{\alpha _{72}} = \frac{2^{72}}{\alpha _{71}} = \frac{2^{71}}{\alpha_{70}} = \frac{2^{\gamma_{_{70}} + 1}}{\alpha_{\gamma_{_{70}}}} = \frac{2,361,183,241,434,822,606,848}{751,587,968,840,192,313,983} \\
&= \underbrace{3.141592653589793238462}_{22\,\,{\rm correct} \,\, {\rm digits}\,\,{\rm of}\,\,\pi}80398052 \ldots
\end{aligned}
\]
remain unchanged. This occurs because the ratio of two adjacent values is 
\[
\alpha_{k + 1} = 2\alpha_k, \qquad 70 \le k \le 74.
\]
However, at $k = 75$ we get
\[
\alpha_{75} = 2\alpha_{74} + 1
\]
since $\alpha_{75} = \alpha_{\gamma_{_{75}}}$ is an odd number. Consequently, with values 
\[
\begin{aligned}
&\alpha_{75} = \alpha_{\gamma_{_{75}}} = 24,050,815,002,886,154,047,457 \\
&\alpha_{76} = 2\alpha_{75} = 48,101,630,005,772,308,094,914 \\
&\alpha_{77} = 4\alpha_{75} = 96,203,260,011,544,616,189,828
\end{aligned}
\]
we can get a slightly more accurate approximation:
\[
\begin{aligned}
\frac{2^{78}}{\alpha_{77}} &= \frac{2^{77}}{\alpha_{76}} = \frac{2^{76}}{\alpha_{75}} = \frac{2^{\gamma_{_{75}} + 1}}{\alpha_{\gamma_{_{75}}}} =\frac{75,557,863,725,914,323,419,136}{24,050,815,002,886,154,047,457} \\
&= \underbrace{3.1415926535897932384626}_{\text{23 correct digits of} \,\, \pi}7335739 \ldots \,\,.
\end{aligned}
\]

These examples showing the relations between the positive integers $\alpha_k$, $\gamma_k$ and $\alpha_{\gamma_{_k}}$ help us to understand the significance of odd integers $\alpha_{\gamma_{_{k}}}$ in the rational approximation \eqref{RA4P} tending to $\pi$ as the integer $k$ increases.

\section{Conclusion}

We present a proof of the irrationality of $\pi$ and show how a contradiction can be reached when the integer $\left\lfloor 2^{k + 1}r\mu_k \right\rfloor - q$ is bounded inside the domain $D$. Examples of the rational approximation tending to $\pi$ with an increasing integer $k$ are provided.

\section*{Acknowledgment}

This work was supported by the National Research Council Canada, Thoth Technology Inc., York University, and Epic College of Technology.



\begin{thebibliography}{999}

\bibitem[Stillwell(2010)]{Stillwell2010}
Stillwell, J. \textit{Mathematics and Its History}, 3rd ed.; Springer: New York, USA, 2010.

\bibitem[Wilson(2018)]{Wilson2018}
Wilson, R. \textit{Euler's Pioneering Equation: The Most Beautiful Theorem in Mathematics}, Oxford University Press: New York, USA, 2018.

\bibitem[Coolidge(1950)]{Coolidge1950}
Coolidge, J.L. \textit{The number e}. {\em Amer. Math. Monthly} {\bf 1950}, {\em 57}(9), 591--602. \url{https://doi.org/10.2307/2308112}.

\bibitem[Davidson(2023)]{Davidson2023}
Davidson, K.R.; Satriano, M. \textit{Integer and Polynomial Algebra}, Mathematical World №31; American Mathematical Society: USA, 2023.

\bibitem[Beckmann(1971)]{Beckmann1971}
Beckmann, P. \textit{A History of Pi}; Golem Press: New York, NY, USA, 1971.

\bibitem[Berggren(2004)]{Berggren2004}
Berggren, L.; Borwein, J.; Borwein, P. \textit{Pi: A Source Book}, 3rd ed.; Springer: New York, NY, USA, 2004.

\bibitem[Agarwal(2013)]{Agarwal2013}
Agarwal, R.P.; Agarwal, H.; Sen, S.K. Birth, growth and computation of pi to ten trillion digits. {\em Adv. Differ. Equ.} {\bf 2013}, {\em 2023}, 100. \url{https://doi.org/10.1186/1687-1847-2013-100}.

\bibitem[Angell(2022)]{Angell2022}
Angell, D. {\textit{Irrationality and Transcendence in Number Theory}}, 1st ed.; CRC Press: Boca Raton, USA, 2022.

\bibitem[Laczkovich(1997)]{Laczkovich1997}
Laczkovich, M. On Lambert's proof of the irrationality of π. {\em Amer. Math. Monthly} {\bf 1997}, {\em 104}(5), 439--443. \url{https://doi.org/10.2307/2974737}.

\bibitem[Zhou(2011)]{Zhou2011}
Zhou, L. Irrationality proofs à la Hermite. {\em Math. Gaz.} {\bf 2011}, {\em 95}(534), 407--413. \url{https://doi.org/10.1017/S0025557200003491}.

\bibitem[Jeffreys(1973)]{Jeffreys2011}
Jeffreys, H. \textit{Scientific Inference}, 3rd ed.; Cambridge University Press: London, UK, 1973.

\bibitem[Niven(1947)]{Niven1947}
Niven, I. A simple proof that π is irrational. {\em Bulletin. Amer. Math. Soc.} {\bf 1947}, {\em 53}(6), 509. \url{https://doi.org/10.1090/s0002-9904-1947-08821-2}.

\bibitem[Huylebrouck(2001)]{Huylebrouck2001}
Huylebrouck, D. Similarities in irrationality proofs for $\pi$, $\ln 2$, $\zeta(2)$, and $\zeta(3)$. {\em Amer. Math. Monthly} {\bf 2021}, {\em 108}(3), 222--231. \url{https://doi.org/10.2307/2695383}.

\bibitem[Bourbaki(2004)]{Bourbaki2004}
Bourbaki, N. {\textit{Functions of a Real Variable: Elementary Theory (Elements of Mathematics)}}, 1st ed.; Springer-Verlag: Berlin Heidelberg, Germany, 2004. \url{https://doi.org/10.1007/978-3-642-59315-4}.

\bibitem[Damini(2020)]{Damini2020}
Damini, D.B.; Dhar, A. How Archimedes showed that π is approximately equal to 22/7. \textit{arxiv} {\bf 2020}, \href{https://arxiv.org/abs/2008.07995}{arXiv:2008.07995}.

\bibitem[Roegel(2020)]{Roegel2020}
Roegel, D. Lambert's proof of the irrationality of pi: context and translation. \textit{HAL open science} \textbf{2020}, \href{https://hal.science/hal-02984214}{hal-02984214}.

\bibitem[Chow(2024)]{Chow2024}
Chow, T.Y. A well-motivated proof that pi is irrational. {\em Hardy-Ramanujan J.} {\bf 2024}, {\em 47}, 26--34. \url{https://doi.org/10.46298/hrj.2025.13361}.

\bibitem[Abrarov(2021)]{Abrarov2021}
Abrarov, S.M.; Siddiqui, R.; Jagpal, R.K.; Quine, B.M. Unconditional applicability of Lehmer's 
measure to the two-term Machin-like formula for pi. {\em Mathematica J.} {\bf 2021}, {\em 23}, 1--23. \url{https://doi.org/doi.org/10.3888/tmj.23–2}.

\bibitem[Abrarov(2022)]{Abrarov2022}
Abrarov, S.M.; Jagpal, R.K.; Siddiqui, R.; Quine, B.M. A new form of the Machin-like formula for $\pi$ by iteration with increasing integers. {{\em J. Integer Seq.} {\bf 2022}, \emph{25}, 22.4.5}. Available online: \url{https://cs.uwaterloo.ca/journals/JIS/VOL25/Abrarov/abrarov5.html} (accessed on 17 March 2026).

\bibitem[Abrarov(2023)]{Abrarov2023}
Abrarov, S.M.; Siddiqui, R.; Jagpal, R.K.; Quine, B.M. A generalized series expansion of the arctangent function based on the enhanced midpoint integration. {\em AppliedMath} {\bf 2023}, {\em 3}, 395--405. \url{https://doi.org/10.3390/appliedmath3020020}.

\bibitem[Abrarov(2018)]{Abrarov2018}
Abrarov, S.M.; Quine, B.M. A formula for pi involving nested radicals. {\em Ramanujan J.} {\bf 2018}, {\em 46}, 657--665. \url{https://doi.org/10.1007/s11139-018-9996-8}.

\bibitem[Servi(2003)]{Servi2003}
Servi, L.D. Nested square roots of 2, {\em Amer. Math. Monthly} {\bf 2003}, {\em 110}(4), 326--330. \url{https://dx.doi.org/10.2307/3647881}.

\bibitem[Abrarov(2025)]{Abrarov2025}
Abrarov, S.M.; Siddiqui, R.; Jagpal, R.K.; Quine, B.M. Application of a new iterative formula for computing π and nested radicals with roots of 2. {\em AppliedMath} {\bf 2025}, {\em 5}(4):156. \url{https://doi.org/10.3390/appliedmath5040156}.

\bibitem[Apostol(2000)]{Apostol2000}
Apostol, T.M. (2000). Irrationality of the square root of two -- a geometric proof. {\em Amer. Math. Monthly}, {\bf 2000}, {\em 107}(9), 841-–842. \url{https://doi.org/10.1080/00029890.2000.12005280}.

\bibitem[Ferreno(2009)]{Ferreno2009}
Ferreño, N.C. Yet another proof of the irrationality of $\sqrt{2}$. {\em Amer. Math. Monthly}, {\bf 2009}, {\em 116}(1), 68--69. \url{https://doi.org/10.1080/00029890.2009.11920911}.

\bibitem[Dmytryshyn(2021)]{Dmytryshyn2021}
Dmytryshyn R.I.; Sharyn S.V. Approximation of functions of several variables by multidimensional S-fractions with independent variables. {\em Carpathian Math. Publ.} {\bf 2021}, {\em 13}(3), 592--607. \url{https://doi.org/10.15330/cmp.13.3.592-607}.

\bibitem[OEIS(1999a)]{OEIS1999a}
The On-Line Encyclopedia of Integer Sequences. OEIS: A024810. Available online: \url{https://oeis.org/A024810} (accessed on 17 March 2026).

\bibitem[Dalzell(1944)]{Dalzell1944}
Dalzell, D.P. On 22/7. {\em J. Lond. Math. Soc.} {\bf 1944}, {\em 19}, 133--134. \url{https://doi.org/10.1112/jlms/19.75_Part_3.133}.

\bibitem[Phillips(1981)]{Phillips1981}
Phillips, G.M.  Archimedes the numerical analyst. {\em Amer. Math. Monthly} {\bf 1981}, {\em 88}(3), 165--169. \url{https://doi.org/10.2307/2320460}.

\bibitem[OEIS(1999b)]{OEIS1999b}
The On-Line Encyclopedia of Integer Sequences. OEIS: A127266. Available online: \url{https://oeis.org/A127266} (accessed on 17 March 2026).

\bibitem[Lorentzen(2008)]{Lorentzen2018}
Lorentzen, L.; Waadeland, H. \textit{Continued Fractions. Volume 1: Convergence Theory}, 2nd edition, Atlantis Studies in Mathematics for Engineering and Science; Atlantis Press: Paris, France, 2008.

\bibitem[Schwerdtfeger(1979)]{Schwerdtfeger1979}
Schwerdtfeger, H. \textit{Geometry of Complex Numbers: Circle Geometry, Moebius Transformation, Non-Euclidean Geometry}; Dover Publications: New York, USA, 1979.

\bibitem[Nield(1982)]{Nield1982}
Nield, D.A. Rational approximations to Pi. {\em New Zealand Math. Mag.} {\bf 1982}, {\em 18}(3), 99--100.

\bibitem[Lucas(2009)]{Lucas2009}
Lucas, S.K. Approximations to $\pi$ derived from integrals with nonnegative integrands. {\em Amer. Math. Monthly} {\bf 2009}, {\em 116}(2), 166--172. \url{https://doi.org/10.1080/00029890.2009.11920923}.

\bibitem[Backhouse(1995)]{Backhouse1995}
Backhouse, N. Pancake functions and approximations to $\pi$. {\em Math. Gaz.} {\bf 1995}, {\em 79}(485), 371--374. \url{https://doi.org/10.2307/3618318}.

\bibitem[Shiu(2014)]{Shiu2014}
Shiu, P. If $e < 2$ then $\pi$ is irrational. {\em Math Gaz.} {\bf 2014}, {\em 98}(543), 497--501. \url{https://doi.org/10.1017/S0025557200008251}.

\bibitem[Rattaggi(2018)]{Rattaggi2018}
Rattaggi, D. Error estimates for the Gregory-Leibniz series and the alternating harmonic series using Dalzell integrals, \href{https://arxiv.org/abs/1809.00998}{{\em arXiv} {\bf 2018}}.

\end{thebibliography}
\end{document}